\documentclass{amsart}

\usepackage{amsmath,xypic,leftidx}
\usepackage{xspace}
\usepackage[psamsfonts]{amssymb}
 
\usepackage[top=1.65in, bottom=1.65in, right=1.3in, left=1.3in]{geometry}

\usepackage{amsfonts}
\usepackage{amssymb}
\usepackage{amsmath}
\usepackage{fancyhdr}
\usepackage{eufrak}
\usepackage{latexsym}
\usepackage{amsbsy}
\usepackage[all]{xy}
\usepackage[latin1]{inputenc}
\usepackage[T1]{fontenc}        
\usepackage{amsmath}

\newtheorem{defn}{Definition}[section]
\newtheorem{thm}[defn]{Theorem}
\newtheorem{prop}[defn]{Proposition}
\newtheorem{lem}[defn]{Lemma}

\newcommand{\fin}{\hfill{\Large$\Box$}\\}

\newcommand{\la}{\lambda}

\newcommand{\C}{\mathbb {C}}

\newcommand{\R}{\mathbb {R}}
\newcommand{\N}{\mathbb {N}}

\newcommand{\Pj}{\mathbb {P}}

\newcommand{\supp}{{\rm supp \, }}

\def\com{\ar@{}[rd]|{\circlearrowleft}}

\DeclareMathOperator{\Ext}{Ext}

\title {Holomorphic motion for Julia sets of holomorphic families of endomorphisms of $\Pj^k$}

\begin{author}[F.~Berteloot]{Fran\c cois Berteloot}
\address{ %
 Universit\'e Paul Sabatier\\
Laboratoire Emile Picard \\
 118, route de Narbonne \\
 31062 Toulouse Cedex \\
  France }
  \email{francois.berteloot@math.univ-toulouse.fr }
\end{author}

\begin{author}[F.~Bianchi]{Fabrizio Bianchi}
\address{ %
 Universit\'e Paul Sabatier\\
Laboratoire Emile Picard \\
 118, route de Narbonne \\
 31062 Toulouse Cedex \\
  France }
  \email{fabrizio.bianchi$@$math.univ-toulouse.fr}
\end{author}

\thanks{This research was partially supported by the ANR project LAMBDA,
ANR-13-BS01-0002. The work of the second author was also partially supported by
the FIRB2012 grant ``Differential Geometry and Geometric Function
Theory''}

\begin{document}

\begin{abstract} 
We build measurable holomorphic motions for Julia sets of holomorphic families of endomorphisms of $\Pj^k$
under various equivalent notions of stability.
\end{abstract}
\maketitle 
\small
{\noindent \emph{Key Words}: holomorphic dynamics, dynamical stability, Lyapunov exponents.

}
 
\normalsize

\section{Introduction and results}

It is a classical result, due to Man\'e-Sad-Sullivan
\cite{MSS_dyn} and Lyubich \cite{L_typical}, that the stability of a holomorphic family of rational maps on the
Riemann sphere is equivalent to the stability of its repelling cycles. This is a fundamental
fact in order to understand the bifurcations within such families. Our aim here is to extend it to
the higher dimensional setting.

A {\it holomorphic family of endomorphisms of} $\Pj^k$ is a holomorphic
map $f : M\times \Pj^k\to M\times \Pj^k$ of the form $f(\la,z)=(\la,f_\la(z))$ where the parameter
space $M$ is a complex manifold and all the maps $f_\la$ have the
same algebraic degree $d$. The {\it Julia set} $J_\la$ of $f_\la$
is, by definition, the support of the maximal entropy measure $\mu_\la$ of
$f_\la$ (see \cite{BD_exp, DS_dyn}). The {\it repelling $J$-cycles of} $f_\la$ are the repelling
cycles of $f_\la$ which belong to $J_\la$. When $k=1$, all repelling cycles are $J$-cycles
but this is not in  general the case when $k>1$ (see \cite{FS_ex}). As it has been shown by
Briend-Duval \cite{BD_exp}, the repelling $J$-cycles of $f_\la$ equidistribute the measure $\mu_\la$.
We will adopt the following definition.

\begin{defn}
The repelling $J$-cycles of $f_\la$  are said to \emph{move holomorphically over}  $M$ if,
for every  $n$, there exists a finite collection of holomorphic maps
$\{\rho_{n,j} : M\to \Pj^k\}_{1\le j\le N_d(n)}$ such that $\{\rho_{n,j} (\la)\}_{1\le j\le N_d(n)}$
is equal to  the set of the $n$-periodic  repelling $J$-cycles of $f_\la$ for each $\la \in M$.
\end{defn}

For  $k=1$, the holomorphic motion of
repelling cycles
yields a holomorphic motion of the Julia sets. This means  that, for a given $\la_0\in M$,
there exists a map $h :M\times J_{\la_0}\to M\times \Pj^1$  of the form $h(\la,z)=(\la,h_\la(z))$
which is continuous, holomorphic in $\la$, one-to-one and which commutes with $f$.
The  map $h$ is the holonomy map of a lamination with transverse measures $\mu_\la$ and can
easily obtained as follows. One first notes that the family $\left(\rho_{n,j} \right)_{j,n}$
is normal (this is actually true even when $k\ge 1$). Then, as the 
$\rho_{n,j} $ have disjoint graphs, Hurwitz lemma implies that sequences like $\gamma_k:=\rho_{j_k,n_k}$
have a unique limit as soon as
$\gamma_k(\la_0)$ converges. It then suffices to set $h(\la,z):=(\la,\lim_k \gamma_k(\la))$
where $z=\lim_k \gamma_k(\la_0)$. This is the so-called $\la$-Lemma.

When $k>1$, the limits of $\gamma_k:=\rho_{j_k,n_k}$ might be not unique and one is
therefore led to consider webs instead of laminations. To this purpose, the following
framework has been introduced in \cite{BD_stab}.
The map $f$ induces a dynamical system $\left({\mathcal J},{\mathcal F}\right)$ where 
\[{\mathcal J}:=\{\gamma : M\to \Pj^k\;\colon\; \textrm{such that}\;f\;\textrm{is holomorphic and}\;
\gamma(\la)\in J_\la\;\textrm{for all}\; \la\in M\}\]
and 
${\mathcal F} : {\mathcal J}\to {\mathcal J}$
is given by${\mathcal F}(\gamma)(\la):=(f_\la\circ \gamma) (\la)$ for all $\la \in M$.
The set $\mathcal J$ is a (possibly empty) metric space for the topology of local uniform convergence. 
Using the compactness properties of probability measures, one easily
sees that if the repelling $J$-cycles of $f$ move holomorphically then there
exists a $\mathcal F$-invariant probability measure $\mathcal M$ on $\mathcal J$ which is compactly
supported and such that $\mu_\la=\int_{\mathcal J} 
\delta_{\gamma(\la)}\; d{\mathcal M}(\gamma)$ for all $\la \in M$. This motivates the following definition.

\begin{defn}
For every $\la \in M$ the projection $p_\la : {\mathcal J}\to \Pj^k$ is
defined by $p_\la(\gamma):=\gamma(\la)$. A  \emph{structural web of} $f$
is a probability measure $\mathcal M$ on $\mathcal J$ such that
${\mathcal F}_\star \left(\mathcal M\right)=\mathcal M $ and $p_{\la\star}\left(\mathcal M\right) =\mu_\la$
for every $\la \in M$.
\end{defn}

It will be crucial in our study to work with  structural webs $\mathcal M$ for which the graphs
$\Gamma_\gamma$ of $\mathcal M$-almost every $\gamma\in {\mathcal J}$ avoid the
grand orbit of the critical set of $f$. The existence of these structural webs
relies on the main results of \cite{BD_stab} and we will show here that some of them are ergodic
(see Proposition \ref{PropGE}).
Using a technique introduced by Briend-Duval in the setting of a single endomorphism of $\Pj^k$ (\cite{BD_exp}),
we will exploit the stochastic properties of $\left({\mathcal J},{\mathcal F},{\mathcal M}\right)$
to show that the iteratated inverse branches of $f$  are exponentially
contracting near the graph
$\Gamma_\gamma:=\{(\la,\gamma(\la))\;\colon\;\la\in M\}$ of $\mathcal M$-almost every
$\gamma\in {\mathcal J}$ (see Theorem \ref{LemSF} and Section \ref{SRIB}).
This implies that for $\mathcal M$-almost every  $\gamma\in {\mathcal J}$
the graph $\Gamma_\gamma$  does not intersect any other graph $\Gamma_{\gamma'}$ where
$\gamma\ne \gamma'\in \supp {\mathcal M}$ (see Proposition \ref{PropExHM}) . 
This approach allows us to
build \emph{measurable holomorphic motions} in the following sense.

\begin{defn}
Let $f:M\times\Pj^k\to M\times\Pj^k$ be a holomorphic family of  endomorphisms of $\Pj^k$ of degree $d\ge 2$.
A \emph{measurable holomorphic motion of the Julia sets} $J_\la$ \emph{over} $M$ is a subset $\mathcal L$ of $\mathcal J$ such 
 ${\mathcal F} \left({\mathcal L} \right)= {\mathcal L}$ and
\begin{itemize}
\item[1)]$\Gamma_\gamma\cap \Gamma_{\gamma'} = \emptyset$ for every distinct $\gamma,\gamma'\in {\mathcal L}$
\item[2)] $\Gamma_\gamma$ does not meet the grand orbit of the critical set of $f$ for every $\gamma\in {\mathcal L}$
\item[3)] $\mu_\la\{\gamma(\la)\;\colon\;\gamma\in {\mathcal L}\} =1$ for every $\la\in M$
\item[4)] the map ${\mathcal F} : {\mathcal L}\to {\mathcal L}$ is $d^k$-to-$1$. 
\end{itemize}
\end{defn}
 Our main result can be stated as follows.
 
\begin{thm}\label{Main}
Let $f:M\times\Pj^k\to M\times\Pj^k$ be a holomorphic family of
endomorphisms of $\Pj^k$ of degree $d\ge 2$,
with $M$ is a simply connected complex manifold. 
If
the $J$-repelling cycles of $f$ move holomorphically over $M$
there exists a measurable holomorphic motion $\mathcal L$
of the Julia sets $J_\la$ over $M$. Moreover, $f$ admits a unique structural
web $\mathcal M$ and ${\mathcal M}\left({\mathcal L}\right)=1$.
\end{thm}

By totally different methods,
Berger and Dujardin (\cite{BD_hyperb}) have recently build measurable
holomorphic motions in the context of polynomial automorphisms of $\mathbb{C}^2$.

When $M$ is a simply connected open subset of the space ${\mathcal H}_d(\Pj^k)$ of 
 endomorphisms of $\Pj^k$ of degree $d \geq 2$, we may
combine this result with some results of \cite{BD_stab} to see that the
holomorphic motion of Julia sets is equivalent to that of repelling cycles.

\begin{thm}\label{main2}  
 Let $f:M\times \Pj^{k}\to M\times \Pj^{k}$
be a holomorphic family of endomorphisms where $M$ is a simply
connected open subset of the space ${\mathcal H}_d(\Pj^k)$ of
 endomorphisms of $\Pj^k$ of degree $d \geq 2$. Let $L(\la)$ denote
 the sum of Lyapunov exponents of $(J_\la,f_\la,\mu_\la)$.
 Then the following assertions are equivalent :
\begin{enumerate}
\item[(A)] the J-repelling cycles move holomorphically
\item[(B)]  the function $L$ is pluriharmonic on $M$
\item[(C)]  there exists a measurable holomorphic motion
of the Julia sets $J_\la$ over $M$.
\end{enumerate}
\end{thm}

\section{From structural webs to holomorphic motions}

In this section, we first show that the existence
of a structural web for which the graphs of almost every $\gamma\in {\mathcal J}$ avoid
the grand orbit of the critical set of $f$ yields an \emph {ergodic} structural web $\mathcal M$ which satisfies the same property.
Then, we explain how the control of iterated inverse branches
in the dynamical system $\left({\mathcal J}, {\mathcal F}, {\mathcal M}\right)$
implies the existence and uniqueness of a holomorphic motion.

In all the section, the critical set of a holomorphic family
$f : M\times \Pj^k\to M\times \Pj^k$  of degree $d$ endomorphisms of $\Pj^k$ will be denoted $C_f$ and 
we set
 $${\mathcal J}_s:=\{\gamma\in {\mathcal J}\;\colon\;
 \Gamma_\gamma\cap \left(\cup_mf^{-m}\left(\cup_n f^n \left(C_f\right)\right) \right) \ne \emptyset\}.$$

\subsection{Existence of good ergodic structural webs}
 
 We aim here to prove the following result.
 
\begin{prop}\label{PropGE} Let $f:M\times\Pj^k\to M\times\Pj^k$ be a holomorphic family of
endomorphisms of $\Pj^k$ of degree $d\ge 2$. If $f$ admits a structural
web ${\mathcal M}_0$ such that ${\mathcal M}_0\left({\mathcal J}_s\right)=0$ then
$f$ admits a structural web ${\mathcal M}'_0$ satisfying ${\mathcal M}'_0\left({\mathcal J}_s\right)=0$
and $\supp {\mathcal M}'_0 \subset \supp {\mathcal M}_0$ and which is ergodic.
\end{prop}

\proof Let us consider the convex set ${\mathcal P}_{web}\left(\mathcal K\right)$ of structural webs of
$f$ which are supported in $\mathcal K$,
where ${\mathcal K}:=\supp \left({\mathcal M}_0\right)$. Note that ${\mathcal F}(\mathcal K)\subset \mathcal K$
since ${\mathcal M}_0$ is $\mathcal F$-invariant. By definition, 
${\mathcal P}_{web}\left(\mathcal K\right)$ is a subset of the set of probability
measures on ${\mathcal K}$. The set  ${\mathcal P}_{web}\left(\mathcal K\right)$ is actually a compact
metric space for the topology of weak convergence of measures. Let us recall why.
Since $\mathcal K$ is a compact subset of the metric space $\mathcal J$, the space $C(\mathcal K)$ of continuous functions
on $\mathcal K$ is a separable Banach space for the norm of uniform convergence. Thus,
by Banach-Alaoglu theorem, the unit ball $B_{C({\mathcal K})'}$ of its dual is metrizable and compact
for the $\textrm{weak}^\star$ topology. Using  the Riesz representation theorem, one then sees
that ${\mathcal P}_{web}\left(\mathcal K\right)$ is closed in $B_{C({\mathcal K})'}$.
Indeed, if ${\mathcal M}_n \in {\mathcal P}_{web}\left(\mathcal K\right)$ converges to ${\mathcal M} \in B_{C({\mathcal K})'}$
then $\mathcal M$ is a probability measure on $\mathcal K$ which is invariant since
${\mathcal F}_\star {\mathcal M}= \lim_n {\mathcal F}_\star {\mathcal M}_n= \lim_n  {\mathcal M}_n= {\mathcal M}$ and is a structural web for $f$ since
$p_{\la \star}{\mathcal M} = \lim_n p_{\la \star}{\mathcal M}_n= \mu_\la$ for every $\la \in M$.

We will use Choquet decomposition theorem to find extremal points ${\mathcal M}'$ in
${\mathcal P}_{web}\left(\mathcal K\right)$ for which ${\mathcal M}'\left({\mathcal J}_s\right)=0$
and then prove the ergodicity of ${\mathcal M}'$ by showing  that these points are also
extremal in the set ${\mathcal P}_{inv}\left(\mathcal K\right)$ of $\mathcal F$-invariant probability measures on $\mathcal K$.

Let us denote by $\Ext\left({\mathcal P}_{web}\left(\mathcal K\right)\right)$ the set of extremal points of the
compact metric space ${\mathcal P}_{web}\left(\mathcal K\right)$.
By Choquet's theorem, there exists a probability
measure $\nu_0$ on $\Ext\left({\mathcal P}_{web}\left(\mathcal K\right)\right)$ such that
\begin{eqnarray*}
{\mathcal M}_0=\int_{\Ext\left({\mathcal P}_{web}\left(\mathcal K\right)\right)} \nu_0\left(\mathcal E\right). 
\end{eqnarray*}
Let us show that the set of structural webs ${\mathcal E}\in \Ext\left({\mathcal P}_{web}\left(\mathcal K\right)\right)$ for which
${\mathcal E}\left({\mathcal J}_s\right) =0$ has full $\nu_0$-measure.
To this purpose we decompose ${\mathcal K}_s:={\mathcal K}\cap{\mathcal J}_s$ into a countable union of compact subsets of $\mathcal J$:
\begin{eqnarray*}
{\mathcal K}_s=\bigcup_n {\mathcal K}_s^n.
\end{eqnarray*}
This can be done by observing that for any relatively compact subset $M'$ of $M$ and
any component of the grand critical orbit of $f$, the set of 
$\gamma \in {\mathcal K}$ for which the graph $\Gamma_\gamma$ intersects that given component over $\overline{M'}$ is compact. 
Now, since ${\mathcal M}_0\left( {\mathcal K}_s^n\right) =0$, it suffices to check that
\begin{eqnarray*}
{\mathcal M}_0\left( {\mathcal K}_s^n\right)=
\int_{\Ext\left({\mathcal P}_{web}\left(\mathcal K\right)\right)} {\mathcal E}\left({\mathcal K}_s^n\right)\;\nu_0\left(\mathcal E\right),
\end{eqnarray*}
which can be seen by approximating the indicatrix $1_{{\mathcal K}_s^n}$ by positive continuous functions on $\mathcal K$.

To conclude the proof we are left to check that any
${\mathcal M}'\in \Ext \left({\mathcal P}_{web}\left(\mathcal K\right)\right)$ is extremal in ${\mathcal P}_{inv}\left(\mathcal K\right)$.
Assume that
${\mathcal M}'=\frac{1}{2} {\mathcal M}_1+\frac{1}{2} {\mathcal M}_2$ where ${\mathcal M}_j\in {\mathcal P}_{inv}\left(\mathcal K\right)$.
Then, as ${\mathcal M}'$ is a structural web for $f$ we have
$\mu_\la=p_{\la\star} \left({\mathcal M}'\right)=
\frac{1}{2} p_{\la\star}\left({\mathcal M}_1\right)+\frac{1}{2} p_{\la\star}\left({\mathcal M}_2\right)$
for every $\la \in M$.
Since $p_\la\circ {\mathcal F}=f_\la\circ p_\la$, the probability measures
$p_{\la\star}\left({\mathcal M}_j\right)$ are $f_\la$-invariant and therefore the ergodicity of $\mu_\la$
implies that $p_{\la\star}\left({\mathcal M}_1\right) =p_{\la\star}\left({\mathcal M}_2\right)=\mu_\la$
for every $\la \in M$. This shows that ${\mathcal M}_1$ and ${\mathcal M}_2$ actually belong to 
${\mathcal P}_{web}\left(\mathcal K\right)$ and the identity ${\mathcal M}'={\mathcal M}_1={\mathcal M}_2$ then
follows from the fact that ${\mathcal M}'$ is extremal in 
${\mathcal P}_{web}\left(\mathcal K\right)$. \fin

\subsection{Iterated inverse branches}\label{SSNE}

In this subsection,  we present the result on the rate of contraction of iterated
inverse branches which will be used to prove the existence of holomorphic motions.
We also fix the framework and some notations for the rest of the paper.

An ergodic structural web $\mathcal M$ for $f$ yields an ergodic dynamical
system $\left({\mathcal J}, {\mathcal F}, {\mathcal M}\right)$,
where $\mathcal F$ is defined by ${\mathcal F}(\gamma)(\la)=f_\la(\gamma(\la))$ for
every $\gamma \in {\mathcal J}$. To study the inverse branches of the map
$\mathcal F$, it is  convenient to transform this system into an injective one.
This is possible using a classical construction called the \emph {natural extension}
which we now describe (we refer to \cite[Chapter 10.4]{CFS_erg} for more details).

Let us assume that $M$ is  simply connected and that ${\mathcal M}\left({\mathcal J}_s\right)=0$.
Recall that ${\mathcal K}:=\supp {\mathcal M}$ is a compact subset of $\mathcal J$.
Setting ${\mathcal X}:={\mathcal K}\setminus{\mathcal J}_s$, it is not difficult to check
that the map ${\mathcal F} : {\mathcal X}\to{\mathcal X}$ is onto.
We may therefore construct the natural extension $\left(\widehat{\mathcal X},\widehat{\mathcal F},
\widehat{\mathcal M}\right)$ of the system $\left({\mathcal X}, {\mathcal F}, {\mathcal M}\right)$ in the following way.
An element of $\widehat{\mathcal X}$ is a left-infinite sequence
$\widehat{\gamma}:=(\cdots,\gamma_{-j},\gamma_{-j+1},\cdots,\gamma_{-1},\gamma_0)$
of elements $\gamma_{j}\in {\mathcal X}$ such that ${\mathcal F}(\gamma_{-j})=\gamma_{-j+1}$
and one defines the map $\widehat{\mathcal F} : \widehat{\mathcal X}\to\widehat{\mathcal X}$ by setting
\[
\begin{aligned}
\widehat{\mathcal F}(\widehat {\gamma})
:&=(\cdots {\mathcal F}(\gamma_j),{\mathcal F}(\gamma_{-j+1}),\cdots,{\mathcal F}(\gamma_{-1}),{\mathcal F}(\gamma_0))\\
& =(\cdots \gamma_{-j+1},\gamma_{-j+2},\cdots, \gamma_0,{\mathcal F}(\gamma_0)).
 \end{aligned}
\]

The map $\widehat{\mathcal F}$ corresponds to the shift operator and is clearly one-to-one and onto.
We may now lift the measure $\mathcal M$
on $\mathcal X$ to a measure $\widehat{\mathcal M}$ on $\widehat{\mathcal X}$ which is characterized by the property that
$$(\pi_{j})_\star\left(\widehat{\mathcal M}\right)={\mathcal M}$$
for any projection $\pi_j:\widehat{\mathcal X}\to\widehat{\mathcal X}$ given by
$$\pi_j(\widehat \gamma)=\pi_j\left((\cdots,\gamma_{-j},\gamma_{-j+1},\cdots,\gamma_{-1},\gamma_0)\right):=\gamma_{-j}.$$
The ergodicity of $\mathcal M$ implies the ergodicity of $\widehat{\mathcal M}$.
We have thus obtained an invertible and ergodic dynamical system
$\left(\widehat{\mathcal X},\widehat{\mathcal F},
\widehat{\mathcal M}\right)$.

For every $\gamma \in {\mathcal J}$ whose graph $\Gamma_\gamma$ does not meet
the critical set of $f$, we denote by $f_\gamma$ the injective map which is
induced by $f$ on some neighbourhood of
 $\Gamma_\gamma$ and by $f^{-1}_{\gamma}$ the inverse branch of $f_\gamma$ which
 is defined on some neighbourhood of $\Gamma_{{\mathcal F}(\gamma)}$. Thus, given
 $\widehat \gamma \in \widehat{\mathcal X}$ and  $n\in \N$ we may define the iterated
 inverse branch $f^{-n}_{\widehat \gamma}$ of $f$ along
 $\widehat \gamma$ and of depth $n$ by
 \[f^{-n}_{\widehat \gamma}:=f^{-1}_{\gamma_{-n}}\circ\cdots \circ f^{-1}_{\gamma_{-2}}\circ f^{-1}_{\gamma_{-1}}.
\]

Let us stress that $ f^{-n}_{\widehat \gamma}$ is defined
on a neighbourhood of $\Gamma_{\gamma_0}$ with values in a 
 neighbourhood of $\Gamma_{\gamma_{-n}}$. Moreover, since only a finite number
 of components of the grand critical orbit of $f$ are involved for 
 defining $ f^{-n}_{\widehat \gamma}$, we may always shrink the parameter space $M$
 to some $\Omega\Subset M$ so  that the domain of definition of  $ f^{-n}_{\widehat \gamma}$
 for a fixed $\widehat \gamma$ 
contains a tubular neighbourhood of $\Gamma_{\gamma_0}\cap \left(\Omega\times \Pj^k\right)$
of the form $T(\gamma_0,\eta) \cap \left(\Omega\times \Pj^k\right)$, where
 $$T(\gamma_0,\eta):=\{(\la,z)\in M\times\Pj^k\;\colon\; d(z,\gamma_0(\la)) <\eta\}.$$

The crucial estimate for our study is given by the following result.

\begin{thm}\label{LemSF}  Let $M$ be a simply connected complex manifold and
$f:M\times\Pj^k\to M\times\Pj^k$ be a holomorphic family of  endomorphisms of $\Pj^k$ of degree $d\ge 2$. Assume that $f$
admits an ergodic structural web $\mathcal M$ such that ${\mathcal M}\left({\mathcal J}_s\right)=0$. Let $\la_0\in M$.
Then there exists a neighbourhood $U_0$ of $\la_0$, an integer $p>0$, a subset $\widehat{\mathcal Y} \subset \widehat{\mathcal X}$
such that $\widehat{\mathcal M}\left(\widehat{\mathcal Y}\right)=1$, a measurable function 
$\widehat{\eta}_p : \widehat{\mathcal Y}\to ]0,1]$ and a constant $A>0$ such that:\\
 for every $\widehat{\gamma}\in \widehat{\mathcal Y}$ and every $n\in p\N^\star$
the iterated inverse branch  $f^{-n}_{\widehat \gamma}$ is defined on the tubular
neighbourhood $T(\gamma_0,\widehat{\eta}_p (\widehat \gamma))\cap (U_0\times\Pj^k)$ of $\Gamma_{\gamma_0}\cap (U_0\times\Pj^k)$
and 
$$f^{-n}_{\widehat \gamma}\left(T(\gamma_0,\widehat{\eta}_p (\widehat \gamma))\cap (U_0\times\Pj^k)\right)
\subset T(\gamma_{-n},e^{-nA})\cap (U_0\times\Pj^k).$$
Moreover, $f^{-n}_{\widehat \gamma}$ is Lipschitz
with $\mbox{Lip} (f^p)_{\widehat{\gamma}}^{-n} \le \widehat{l}_p (\widehat{\gamma}) e^{-nA}$,
where $\widehat{l}^p(\widehat{\gamma})~\ge~1$.
\end{thm}

The last section of the paper will be devoted to the proof of the above result.

\subsection{Existence and uniqueness of holomorphic motions}

The existence of a unique holomorphic motion will basically be obtained by combining Proposition \ref{PropGE} with the following one.

\begin{prop}\label{PropExHM}
Let $M$ be a simply connected complex  manifold and $f:M\times\Pj^k\to M\times\Pj^k$ be
a holomorphic family of  endomorphisms of $\Pj^k$ of degree $d\ge 2$.
Assume that there exists an ergodic structural web ${\mathcal M}_0$  of $f$
such that ${\mathcal M}_0\left({\mathcal J}_s\right)=0$ and let ${\mathcal K}_0:=\supp {\mathcal M}_0$.
 Then
$$
{\mathcal M}_0\left(\{\gamma\in {\mathcal K}_0\;\colon\; \exists k\in \N,
\exists \gamma'\in {\mathcal K}_0\;\textrm{s.t.}\;\Gamma_{{\mathcal F}^k(\gamma)} \cap \Gamma_{\gamma'} \ne \emptyset\;
\textrm{and}\; {\mathcal F}^k(\gamma)\ne \gamma'\}\right)=0.
$$
Moreover, any other structural web of $f$ coincides with ${\mathcal M}_0$ as a probability measure on ${\mathcal J}$.
\end{prop}

\proof to prove the first statement, it
is sufficent  to show that for any fixed $k\in \N$  and any $\la_0 \in M$
there exists a neighbourhood $U_0$ of $\la_0$ such that
\begin{eqnarray}\label{firstgoal}
{\mathcal M}_0\left(\{\gamma\in {\mathcal K}_0\;\colon\; \exists \gamma'\in {\mathcal K}_0\;
\textrm{s.t.}\;\Gamma_{{\mathcal F}^k(\gamma)} \cap \Gamma_{\gamma'}\cap \left(U_0\times\Pj^k\right) \ne \emptyset\;
\textrm{and}\; {{\mathcal F}^k(\gamma)} \ne \gamma'\}\right)=0
\end{eqnarray}

 To this purpose, we shall work with the natural extension $\left(\widehat{\mathcal X},\widehat{\mathcal F},
\widehat{{\mathcal M}_0}\right)$ of the system $\left({\mathcal X}, {\mathcal F}, {\mathcal M}_0\right)$
which has been constructed in subsection \ref{SSNE} and apply
Theorem \ref{LemSF}. Let $U_0$ be a neighbourhood of $\la_0$ given by that theorem;
we may assume that $U_0$ is simply connected and that $U_0 \Subset M$.
We recall that ${\mathcal X}\subset {\mathcal K}_0$ and set ${\mathcal F}^k(\gamma)=:\gamma_k$
and $\widehat{\mathcal F}^k(\widehat{\gamma})=:\widehat{\gamma}_k$.

For any $B\subset U_0$, we define the ramification functions $R_B$ by setting
\begin{eqnarray*}
R_B(\gamma):=\sup_{\gamma'\in  {\mathcal K}_0 : \Gamma_{\gamma'\vert B} \cap \Gamma_{\gamma\vert B} \neq \emptyset}
\sup_B d\left(\gamma(\la),\gamma'(\la)\right), \;\;
 \forall \gamma\in {\mathcal J}.
\end{eqnarray*}
Let $p$ be the integer and $\widehat{\eta}_p:\widehat{\mathcal Y} \to ]0,1]$ be
the measurable function given by  Theorem \ref{LemSF}. For $\epsilon >0$ we set
 \begin{eqnarray*}
\widehat{\mathcal Y}_\epsilon:=\{\widehat \gamma\in \widehat{\mathcal Y}\;\colon\;
\widehat{\eta}_p({\widehat \gamma}_k) >\epsilon\;\textrm{and}\; R_{U_0}(\gamma_k)>0\}.
\end{eqnarray*}

It then  suffices to prove that $\widehat{{\mathcal M}_0}\left( \widehat{\mathcal Y}_\epsilon\right)=0$
for every $\epsilon >0$ as it follows from the following observation:
 \begin{eqnarray*}
 {\mathcal M}_0\left(\{\gamma\in {\mathcal K}_0\;\colon\;
 \exists \gamma'\in {\mathcal K}_0\;\textrm{s.t.}\;\Gamma_{\gamma_k}
 \cap \Gamma_{\gamma'}\cap\left(U_0\times\Pj^k\right) \ne \emptyset\;
\textrm{and}\; \gamma_k\ne \gamma'\}\right)\\
={\mathcal M}_0\left(\{\gamma\in{\mathcal K}_0\;\colon\; R_{U_0}(\gamma_k)>0\}\right)
 = {\mathcal M}_0\left(\{\gamma\in{\mathcal X}\;\colon\; R_{U_0}(\gamma_k)>0\}\right)
\\
=\widehat{{\mathcal M}_0}\left(\{\widehat{\gamma}\in{\widehat {\mathcal Y}}\;\colon\; R_{U_0}(\gamma_k)>0\}\right)=
\widehat{{\mathcal M}_0}\left(\cup_{\epsilon >0} \widehat{\mathcal Y}_\epsilon\right).
\end{eqnarray*}

Let us proceed by contradiction and assume that
$\widehat{{\mathcal M}_0}\left( \widehat{\mathcal Y}_\epsilon\right) >0$ for some $\epsilon>0$.
Owing to the equicontinuity of ${\mathcal X}$ (we recall that ${\mathcal X} \subset \supp {\mathcal M}_0$)
we may cover $U_0$ with  finitely many open sets $B_i \subset U_0$, say with $1\le i\le N$, such that

\begin{eqnarray}\label{equic}
\forall \gamma,\gamma'\in {\mathcal X},\;\forall \la_1\in B_i\;:\; \gamma(\la_1)=\gamma'(\la_1)
\Rightarrow \sup_{\la\in B_i} d\left(\gamma(\la),\gamma'(\la)\right) < \epsilon.
\end{eqnarray}

As $R_{U_0}(\gamma)=0$ when $\max_{1\le i\le N} R_{B_i}(\gamma) =0$
(by analyticity we have $\gamma=\gamma'$ on $U_0$ if 
$\gamma=\gamma'$ on some $B_i$), there exists $1\le j\le N$ and $\alpha>0$ such that:

\begin{eqnarray*}
\widehat{{\mathcal M}_0}\left(\{ {\widehat \gamma}\in \widehat{\mathcal Y}\;\colon\;
\widehat{\eta}_p({\widehat \gamma}_k) >\epsilon\;\textrm{and}\;R_{B_j}(\gamma_k)>\alpha\}\right) >0.
\end{eqnarray*}

Let us set $\widehat{\mathcal Y}_{\epsilon,j,\alpha}:=\{\widehat \gamma\in \widehat{\mathcal Y}\;\colon\;
\widehat{\eta}_p({\widehat \gamma}_k) >\epsilon\;\textrm{and}\; R_{B_j}(\gamma_k)>\alpha\}.$
Applying Poincar\'e recurrence theorem to $\widehat{\mathcal F}^{-p}$, we
find $\widehat\gamma\in \widehat{\mathcal Y}_{\epsilon,j,\alpha}$ and an increasing
sequence of integers $(n_q)_q$ with $n_q\in p\N$ such that
$\widehat{\mathcal F}^{-n_q} (\widehat \gamma) \in \widehat{\mathcal Y}_{\epsilon,j,\alpha}$
fo every $q\in \N$. In particular $\widehat\gamma\in \widehat{\mathcal Y}_{\epsilon}$
and $R_{B_j} (\gamma_{k-n_q})>\alpha$ for every $q\in \N$.
We will reach a contradiction by establishing that
\begin{eqnarray}\label{Lim0}
\lim_{m\to+\infty} R_{B_i} (\gamma_{k-m p})=0, \;\;\forall i\in \{1,\cdots,N\},\;
\forall \widehat\gamma\in \widehat{\mathcal Y}_{\epsilon}.
\end{eqnarray}

To this purpose we shall use Theorem \ref{LemSF} to show that
$R_{B_i} (\gamma_{k-n})\le e^{-nA}$ when $n\in p\N$ and
$\widehat\gamma\in \widehat{\mathcal Y}_{\epsilon}$. Let $\gamma'~\in~{\mathcal K}_0$ such
that $\gamma'(\la_1)=\gamma_{k-n} (\la_1)$ for some $\la_1\in B_i$. Then
$({\mathcal F}^n\gamma')(\la_1)=\gamma_k(\la_1)$ and thus, according to (\ref{equic}), 
$\sup_{\la\in B_i} d\left(({\mathcal F}^n\gamma')(\la),\gamma_k(\la)\right)
< \epsilon < \widehat{\eta}_p({\widehat \gamma}_k)$. This means that
\begin{eqnarray}\label{tube}
\Gamma_{{\mathcal F}^n\gamma'} \cap \left(B_i\times \Pj^k\right)\subset
T\left(\gamma_k,\widehat{\eta}_p({\widehat \gamma}_k)\right)\cap \left(B_i\times \Pj^k\right).
\end{eqnarray}

Now, according to Theorem \ref{LemSF}, the inverse branch $f_{{\widehat \gamma}_k}^{-n}$ of $f^n$
is defined on
$T\left(\gamma_k,\widehat{\eta}_p({\widehat \gamma}_k)\right)\cap \left(U_0\times \Pj^k\right)$
and maps it biholomorphically into the tube $T\left(\gamma_{k-n},e^{-nA}\right)$.
As $B_i \subset U_0$, this yields:
\begin{eqnarray}\label{tubeagain}
f_{{\widehat \gamma}_k}^{-n}\left(
T\left(\gamma_k,\widehat{\eta}_p({\widehat \gamma}_k)\right)\cap \left(B_i\times \Pj^k\right) \right)
\subset T\left(\gamma_{k-n},e^{-nA}\right).
\end{eqnarray}

By construction we have $f_{{\widehat \gamma}_k}^{-n}\left(
\Gamma_{\gamma_k}\right) =\Gamma_{\gamma_{k-n}}$ and thus 
$f_{{\widehat \gamma}_k}^{-n}\left(({\mathcal F}^n\gamma')(\la_1)\right)
=f_{{\widehat \gamma}_k}^{-n}\left(\gamma_k(\la_1)\right) =\gamma_{k-n}(\la_1)=\gamma'(\la_1)$.
This implies that $f_{{\widehat \gamma}_k}^{-n}\left(\Gamma_{{\mathcal F}^n\gamma'}\right)=\Gamma_{\gamma'}$
which in turns,
by 
(\ref{tube}) and (\ref{tubeagain}),
implies that,
$\sup_{\la\in B_i} d(\gamma'(\la),\gamma_{k-n}(\la) \le e^{-nA}$.
Then (\ref{Lim0}) follows and  (\ref{firstgoal}) is proved.

We now prove the uniqueness assertion. Let us fix $\la \in M$ and,
for any Borel subset $\mathcal A$ of $\mathcal J$, let us set
${\mathcal A}_\la:=\{\gamma(\la)\;\colon\;\gamma\in {\mathcal A}\}$. Then, as
${\mathcal A}\subset p_\la^{-1}\left({\mathcal A}_\la\right)$ we have
$$\mu_\la({\mathcal A}_\la) = (p_{\la\star} {\mathcal M})({\mathcal A}_\la) =
{\mathcal M}\left( p_\la^{-1}\left({\mathcal A}_\la\right)\right) \ge {\mathcal M}({\mathcal A})$$
for every structural web $\mathcal M$ of $f$.
On the other hand, it follows from (\ref{firstgoal}) applied for $k=0$ that 
$$\mu_\la({\mathcal A}_\la) ={\mathcal M}_0\left( p_\la^{-1}\left({\mathcal A}_\la\right)\right) = {\mathcal M}_0({\mathcal A}).$$
 We thus have
${\mathcal M}_0({\mathcal A}) \ge {\mathcal M}({\mathcal A})$ for any borelian subset
$\mathcal A$ of $\mathcal J$ and this implies that the measures $\mathcal M$ and ${\mathcal M}_0$
must coincide since  both are probability measures on 
$\mathcal J$.
 \fin

We now give a proof of Theorem \ref{Main}.

\proof  By assumption,
for every  $n\ge 1$ we have subsets ${\mathcal R}_n:=\{\rho_{n,j} \;\colon \; 1\le j\le N_d(n)\}$ of $\mathcal J$
such that the 
$\rho_{n,j}(\la)$ are
the
repelling $n$-periodic points of $f_\la$ for every $\la\in M$.
Notice that $N_d(n)\le d^{kn}$ since any $f_\la^n$ cannot have more than $d^{kn}$ fixed points, and moreover that
$\lim_n d^{-kn}N_d(n) =1$.
Without loss of generality we may assume that ${\mathcal F}\left({\mathcal R}_n\right) ={\mathcal R}_n$.
We define a sequence $\left({\mathcal M}_n\right)_n$ of 
$\mathcal F$-invariant discrete probability measures on $\mathcal J$ by setting
${\mathcal M}_n:=\frac{1}{N_d(n)} \sum_{j=1}^{N_d(n)} \delta_{\rho_{n,j}(\la) }$.

Considering lifts to $\C^{k+1}$, one may check that the family $\cup_n\{\rho_{n,j}\;\colon \; 1\le j\le N_d(n)\}$
is equicontinuous.
Thus, there exists a compact subset $\mathcal K$ of $\mathcal J$ such that
${\mathcal F}\left(\mathcal K\right)\subset {\mathcal K}$ and $\supp M_n \subset {\mathcal K}$ for every $n\ge 1$.
Therefore, after taking a subsequence, 
$\left({\mathcal M}_n\right)_n$  converges to a compactly supported probability measure $\mathcal M$ on $\mathcal J$.

Let us check  that $\mathcal M$ is a structural web for $f$. For every $\la \in M$, we have 
 $p_{\la\star} \left({\mathcal M}\right)= p_{\la\star}\left(\lim_n {\mathcal M}_n\right) = \lim_n p_{\la\star}\left( 
{\mathcal M}_n\right)=
\lim_n \frac{1}{N_d(n)} \sum_{j=1}^{N_d(n)} \delta_{\rho_{n,j}(\la) }= \mu_\la$
where the last equality follows from $\lim_n d^{-kn}N_d(n) =1$ and Briend-Duval
equidistribution theorem \cite{BD_exp}.
The invariance of $\mathcal M$ follows from that of the ${\mathcal M}_n$.
Indeed,
\[{\mathcal F}_\star \left(\mathcal M\right) = {\mathcal F}_\star \left(\lim_n {\mathcal M}_n\right)
= \lim_n {\mathcal F}_\star \left( {\mathcal M}_n\right)\\ =  \lim_n  \left( {\mathcal M}_n\right) ={\mathcal M}.\]

We shall now prove that ${\mathcal M}\left({\mathcal J}_s\right)=0$. To this purpose
we first notice that, according to \cite[Proposition 3.1]{BD_stab},
$dd^c L(\la)=0$ where
$L(\la)$ is the sum of Lyapunov exponents of $f_\la$ with respect to the measure
$\mu_\la$. Then we use Proposition 3.3 of \cite{BD_stab} to deduce  that $M$ does not contain
Misiurewicz parameters. recall that a parameter $\la_0$ is called Misiurewicz if
there exists a holomorphic map $\rho : N_{\la_0}\to \Pj^k$ defined on some neighbourhood
$N_{\la_0}$ of $\la_0$ such that
$\rho(\la)$ is a repelling $J$-cycle of $f_\la$ for every $\la\in N_{\la_0}$, the
point $(\la_0,\rho(\la_0))$ belongs to $f^k(C_f)$ for some $k\ge 1$ but 
the graph of $\rho$ is not contained in $f^k(C_f)$.

We can now see  that for every $k \in \N$ and every $\gamma\in \supp \mathcal M$ one has:
$$\Gamma_\gamma \cap f^k(C_f) \ne \emptyset \Rightarrow \Gamma_\gamma\subset f^k(C_f).$$
 Indeed, if this were not  the case, by Hurwitz theorem, we could find
 some $\gamma' \in \cup_n \supp {\mathcal M}_n$ such that $\Gamma_\gamma \cap f^k(C_f)
 \ne \emptyset$ and $\Gamma_\gamma$ is not contained in $f^k(C_f)$.
When $k=0$ this is clearly impossible since $\gamma'(\la)$ is a repelling cycle
of $f_\la$ and when $k\ge 1$, this is impossible because $M$ does not contain
Misiuerewicz parameter.

So, fixing any $\la_0\in M$ we get 
\begin{eqnarray*}
\begin{aligned}
{\mathcal M} \left(\{ \gamma\in {\mathcal J}\;\colon\; \Gamma_\gamma\cap \left(\cup_{k\ge 0} f^k(C_f)\right) \ne \emptyset\}\right)
& = {\mathcal M} \left(\{ \gamma\in {\mathcal J}\;\colon\; \Gamma_\gamma\subset \left(\cup_{k\ge 0} f^k(C_f)\right) \}\right) \\
 & \le {\mathcal M} \left(\{ \gamma\in {\mathcal J}\;\colon\; (\la_0,\gamma(\la_0)) \in 
 \left(\cup_{k\ge 0} f^k(C_f)\right) \}\right)\\
 & = \mu_{\la_0}  \left(\cup_{k\ge 0} f_{\la_0}^k(C_{f_{\la_0}})\right) =0
\end{aligned}
\end{eqnarray*}
where the two last equalities come from $p_{\la_0\star}\left({\mathcal M}\right)=\mu_{\la_0}$
and the fact that $\mu_{\la_0}$ does not charge pluripolar sets in $\Pj^k$. The estimate ${\mathcal M}\left({\mathcal J}_s\right)=0$
finally follows from the $\mathcal F$-invariance of $\mathcal M$.

According to Proposition \ref{PropGE}, there exists a structural
web ${\mathcal M}_0$ which is ergodic as $\mathcal F$-invariant measure on ${\mathcal K}_0:=\supp {\mathcal M}_0$ and such that
${\mathcal M}_0\left({\mathcal J}_s\right)=0$. We 
consider $${\mathcal L}_+:=\{\gamma\in {\mathcal K}_0\setminus {\mathcal J}_s\;\colon\; \forall \gamma'\in {\mathcal K}_0, \forall k \in \N,  
\Gamma_{{\mathcal F}^k(\gamma)}\cap \Gamma_{\gamma'} \ne \emptyset \Rightarrow {\mathcal F}^k(\gamma)=\gamma'\}.$$
By Proposition \ref{PropExHM} we have ${\mathcal M}_0\left({\mathcal L}_+\right)=1$ and, by construction,
${\mathcal L}_+$ satisfies the following 
properties:
\begin{itemize}
\item[1)] ${\mathcal L}_+\subset {\mathcal J}\setminus {\mathcal J}_s$
\item[2)] ${\mathcal F}\left({\mathcal L}_+\right)\subset{\mathcal L}_+$
\item[3)] $ \forall \gamma, \gamma'\in {\mathcal L}_+\; : \; \Gamma_\gamma \cap \Gamma_{\gamma'}
\ne \emptyset \Rightarrow \gamma=\gamma'$.
\end{itemize}
The set ${\mathcal L}:=\cup_{m\ge 0} {\mathcal F}^{-m} \left({\mathcal L}_+\right)$
satisfies the same properties and moreover ${\mathcal F} : {\mathcal L}\to {\mathcal L}$ is $d^k$-to-1. \fin

\section{Rate of contraction of iterated inverse branches}\label{SRIB}
The aim of this section is to prove Theorem \ref{LemSF}.
Since everything here is local, we may assume that the parameter space $M$ is
an open subset of $\C^m$ which we endow with the euclidean norm.

In all the section, we shall use the notations and concepts introduced in
Subsection \ref{SSNE} and in particular work with the natural
extension $\left(\widehat {\mathcal X},\widehat{\mathcal F}, \widehat{\mathcal M}\right)$
where $\mathcal M$ is an ergodic structural web such that ${\mathcal M}({\mathcal J}_s)=0$ and 
 ${\mathcal X}=\supp {\mathcal M} \setminus {\mathcal J}_s={\mathcal K}\setminus {\mathcal J}_s$.
 Let us recall that, by definition of a structural web, the support $\mathcal K$ of $\mathcal M$
 is a compact subset of $\mathcal J$.

\subsection{Lyapunov exponent and rate of contraction}

The material presented in this subsection is not original.
We simply adapt to the context of the system $({\mathcal J},{\mathcal F},{\mathcal M})$ the
tools which have been first introduced in \cite{BD_exp} by Briend-Duval for the case of
a single holomorphic endomorphism of $\Pj^k$.

We first need to fix good sets of holomorphic charts on $\Pj^k$. 
For any $\tau >0$, we may find a covering $\Pj^k=\cup_{i=1}^N V_i$ by open sets and a collection of holomorphic maps 
$$\psi_i : V_i \times B_{\C^k} (0,R_0) \to \Pj^k$$
such that $\psi_{i,x}:=\psi_i(x,\cdot)$ is a (holomorphic)
chart centered at $x\in V_i$ (which means that  $\psi_{i,x}: B_{\C^k}(0,R_0)\to \psi_{i,x}\left(  B_{\C^k}(0,R_0) \right)$
is a biholomorphism and $\psi_{i,x}(0)=x$) and 
\begin{eqnarray}\label{LipChart}
e^{-\tau /2} \vert z-z'\vert \le d_{\Pj^k}\left( \psi_{i,x}(z), \psi_{i,x}(z')\right) \le e^{\tau/2} \vert z-z'\vert
\end{eqnarray}
for  all $(x,z)\in V_i \times B_{\C^k}(0,R_0)$ and all $1\le i\le N$.

We will now use these holomorphic charts to express the
restrictions of $f^n$ on suitable neighbourhoods of graphs $\Gamma_\gamma$.
Let us fix $\la_0$ in $M$. Since the family $\mathcal K$ is locally equicontinuous,
there exists a relatively compact open ball $W_0$ centered at  $\la_0$ in $M$  such that:
$$\forall \gamma \in {\mathcal K},\;\exists i\in\{1,2,\cdots,N\}\;\textrm{such that}\;
\gamma(\la)\in V_i \;\textrm{for all}\; \la \in  \overline {W_0}.$$
For all $\gamma \in {\mathcal K}$ we set 
$$i(\gamma):= \inf \{1\le i\le N\;\colon\; \gamma(\la)\in V_i\;\textrm{for all}\; \la \in  \overline{W_0}\}.$$
Then, for every $n\ge 1$ there exists  $R_n\in ]0,R_0]$ such
that the maps $F_{\gamma(\la)}^n(z)$ given by
$$F_{\gamma(\la)}^n(z):= \left(\psi_{i({\mathcal F}^n(\gamma)),{\mathcal F}^n(\gamma)(\la)} \right)^{-1}
\circ f_{\la}^n \circ \psi_{i(\gamma),\gamma(\la)} (z)$$
are well defined and holomorphic on $\overline{W_0}\times \overline{B_{\C^k}(0,R_n)}$
(i.e. on a fixed neighbourhood of $\overline{W_0}\times \overline{B_{\C^k}(0,R_n)}$) for every $\gamma \in {\mathcal K}$.
This follows immediately from the definitions and the uniform continuity of $f^n$ on $\overline{W_0}\times \Pj^k$.

 As 
$F_{\gamma(\la)}^n(z)$ is locally invertible at the origin when $\gamma\notin {\mathcal J}_s$, we may now define
functions $u_n$ on ${\mathcal X}\times {W_0}$ by setting
$$u_n(\gamma,\la):=\ln \Vert  (DF_{\gamma(\la)}^n(0))^{-1}\Vert;\;\;\forall (\gamma,\la)
\in {\mathcal X}\times \overline{W_0},\;\forall n\ge 1.$$
Let us stress that $ (DF_{\gamma(\la)}^n(0))^{-1}$ depends holomorphically on $\la \in {W_0}$.

From now on we consider  three open balls $U_0\Subset V_0\Subset W_0$ centered at $\la_0$ in $M$.
Let us  introduce the function $r_p$  on $\mathcal X$ by setting:
$$ r_p(\gamma) :=e^{-2\sup_{\la \in U_0} u_p(\gamma,\la)}.$$ 

The next lemma shows that the function $r_p$  measures the size of tubular neighbourhoods of $\Gamma_\gamma$
on which $f^p$ is invertible and also the rate of contraction of the inverse branch $(f^p)_{\gamma}^{-1}$. 

\begin{lem}\label{Lemr} Let $r_p : {\mathcal X} \to \R^{+\star}$ be the function
defined by  $r_p(\gamma) := e^{-2\sup_{\la \in {U_0}} u_p(\gamma,\la)}$.
For any sufficently small  $\epsilon >0$ there exists $C_p(\epsilon) >0$ such that
for any $\gamma\in {\mathcal X}$ the map $f^p$ admits an inverse branch
$(f^p)_{\gamma}^{-1}$ on the tube $ T({\mathcal F}^p(\gamma),C_p(\epsilon)r_p(\gamma))\cap(U_0\times\Pj^k)$ which maps 
$\Gamma_{{\mathcal F}^p(\gamma)}\cap(U_0\times\Pj^k)$ to $\Gamma_\gamma\cap(U_0\times\Pj^k)$ and is Lipschitz with
$\mbox{Lip}\left((f^p)_{\gamma}^{-1}\right) \le e^{\tau + \epsilon/3} r_p(\gamma)^{-1/2}$ .
\end{lem}

\proof We apply  a quantitative version of the inverse mapping theorem to the maps
 $F_{\gamma(\la)}^p$ where  $\gamma\in {\mathcal X}$ and $\la\in {U_0}$. Our
 precise references for that are \cite[Lemma 2]{BD_exp} and \cite[Lemma 1.1.32]{D_propr}.
 This yields $\delta_p(\epsilon)>0$ such that the following statement holds for any $(\gamma,\la) \in {\mathcal X}\times {U_0}$:
 
\begin{center} the inverse  $\left(F_{\gamma(\la)}^p\right)^{-1}$ of $F_{\gamma(\la)}^p$
is defined on $B_{\C^k} \left(0,\delta_p(\epsilon)\Vert (DF_{\gamma(\la)}^p(0))^{-1}\Vert^{-2} \right)$
and is Lipschitz with
$\mbox{Lip}\left(F_{\gamma(\la)}^p\right)^{-1} \le e^{\frac{\epsilon}{3}} \Vert (DF_{\gamma(\la)}^p(0))^{-1}\Vert$.
\end{center}

Let us stress that  $\delta_p(\epsilon)$ does not depend on $(\gamma,\la)$ since
one may take $\delta_p(\epsilon)=\frac{R_p(1-e^{-\frac{\epsilon}{3})}}{M}$ where
$M=\sup_{\la\in {U_0},\gamma\in{\mathcal X}} \Vert F_{\gamma(\la)}^p\Vert_{{\mathcal C}^2,\overline{B(0,R_p)}}$.

Taking into account the fact that the distorsions introduced by the charts are bounded
(see (\ref{LipChart})), one sees that a similar statement holds for $f^p_\la$ after dividing
$\delta_p(\epsilon)$ by a suitable large constant and modifying the Lipschitz estimate. Namely:
\begin{center}
 for every $\gamma\in{\mathcal X}$ and every $\la\in {U_0}$, $f^p_\la$ admits an
 inverse branch $(f^p)_{\gamma(\la)}^{-1}$ on
 $B_{\Pj^k} \left(f^p_\la(\gamma(\la)),\delta_p(\epsilon)\Vert (DF_{\gamma(\la)}^p(0))^{-1}\Vert^{-2} \right)$
 which maps $f^p_\la(\gamma(\la))$ to
$\gamma(\la)$ and is Lipschitz with
$\mbox{Lip} (f^p)_{\gamma(\la}^{-1} \le e^{\tau + \frac{\epsilon}{3}} \Vert (DF_{\gamma(\la)}^p(0))^{-1}\Vert$.
\end{center}
As $f^p(\la,z)=(\la,f^p_\la(z))$ is a fibered map, the statement
follows since 
$\inf_{\la\in {U_0}} \Vert (DF_{\gamma(\la)}^p(0))^{-1}\Vert^{-2}=e^{-2\sup_{\la \in {U_0}} u_p(\gamma,\la)}=r_p(\gamma)$.
\fin

We  now define  functions $\widehat{u}_n$ on $\widehat{\mathcal X}$ by setting:
\begin{eqnarray}\label{new}
\widehat{u}_n(\widehat{\gamma}):=\sup_{\la\in U_0} u_n(\gamma_0,\la);\;\forall \widehat{\gamma}\in \widehat{\mathcal X}.
\end{eqnarray}

The stochastic properties of the functions $u_n$ actually allow to control
the asymptotic behaviour of $r_p(\gamma_{-n})$ for almost every
$\widehat{\gamma}$ in $\widehat {\mathcal X}$ and to obtain the estimate given in
Theorem \ref{LemSF}. This is what our next Proposition shows.

\begin{prop}\label{PropInter}
 Let $f:M\times\Pj^k\to M\times\Pj^k$ be a holomorphic family of
 endomorphisms of $\Pj^k$ of degree $d\ge 2$ which  admits an ergodic
 structural web $\mathcal M$ such that ${\mathcal M}\left({\mathcal J}_s\right)=0$.
 
 Let $\tau >0$ and $\epsilon>0$ be such that Lemma \ref{Lemr} holds and $-\frac{\ln d}{2} +\tau + 2\epsilon <0$.
  Let $U_0\Subset V_0\Subset W_0$ be relatively compact open balls centered
  at $\la_0 \in M$ and let 
  $\widehat{u}_n$ be the functions defined on  $\widehat{\mathcal X}$ as in (\ref{new}).
  
 Assume that the functions $\widehat{u}_n$ are $\widehat{\mathcal M}$-integrable and that there exists $L\le \frac{-\ln d}{2}$ such that 
 $$\lim_n\frac{1}{n} \int_{\widehat{X}} \widehat{u}_n\;d\widehat{\mathcal M}=L.$$

 Then there exists an integer $p$, a subset $\widehat{\mathcal Y} \subset \widehat{\mathcal X}$
such that $\widehat{\mathcal M}\left(\widehat{\mathcal Y}\right)=1$, a measurable function 
$\widehat{\eta}_p: \widehat{\mathcal Y}\to ]0,1]$ and a constant $A>0$ such that:

for every $\widehat{\gamma}\in \widehat{\mathcal Y}$ and $n\in p\N^\star$
the iterated inverse branch  $f^{-n}_{\widehat \gamma}$ is defined on
the tubular neighbourhood $T(\gamma_0,\widehat{\eta}_p (\widehat \gamma))\cap (U_0\times\Pj^k)$
of $\Gamma_{\gamma_0}\cap (U_0\times\Pj^k)$ and 
$f^{-n}_{\widehat \gamma}\left(T(\gamma_0,\widehat{\eta}_p (\widehat \gamma))\cap (U_0\times\Pj^k)\right)
\subset T(\gamma_{-n},e^{-nA})\cap (U_0\times\Pj^k)$.
\end{prop}

 \proof  We shall use the following Lemma for which we refer to \cite{D_propr}.
 
  \begin{lem}\label{lemDD}
 Let $\varphi$ and $\psi_n$ be measurable and strictly positive functions on a
 metric space $X$. Suppose that $\lim_n \frac{1}{n} \ln \psi_n(x)=0$ for all $x\in X$. Then
 for any $\epsilon >0$ there exists measurable functions $\alpha,\beta : X\to ]0,+\infty[$ such that
 $\alpha e^{-n\epsilon} \le \psi_n\le \beta e^{n\epsilon}$ for all $n\in \N$  and $\alpha\le \varphi \le \beta$.
 \end{lem}
 
 As by assumption  $\lim_n \frac{1}{n}\int_{\widehat{\mathcal X}} \widehat{u_n}\;d\widehat{\mathcal M} =L$
 with $L\le -\frac{\ln d}{2}$,  we may pick $p\in \N^\star$ such
 that $\frac{1}{p}\int_{\widehat{\mathcal X}} \widehat{u}_p\;d\widehat{\mathcal M} =:L''\le L +\epsilon$.
 Setting $L':=pL''$, we note that $L'+\tau + \epsilon <0$.
 
 Then we can apply Birkhoff ergodic theorem to find  $\widehat{\mathcal Y} \subset\widehat{\mathcal X}$
 such that $\widehat{\mathcal M} (\widehat{\mathcal Y})=1$ and
 \begin{eqnarray}\label{Bir}
 \lim_n \frac{1}{n} \sum_{j=1}^n \widehat{u}_p\left(\widehat{\mathcal F}^{-j} (\widehat{\gamma})\right) 
= \int_{\widehat{\mathcal X}} \widehat{u}_p\;d\widehat{\mathcal M} =pL''=L';\;\;\forall \widehat{\gamma}\in \widehat{\mathcal Y}.
  \end{eqnarray}

 By (\ref{Bir}) we have
 $\lim_n \frac{1}{n} \ln r_p(\gamma_{-n}) =\lim_n \left( -2\widehat{u}_p(\widehat{\mathcal F}^{-n}  (\widehat{\gamma}) )\right)=0$
 for all $\widehat{\gamma}\in\widehat{\mathcal Y}$. Thus, applying Lemma \ref{lemDD}
 to the functions $\psi_n(\widehat{\gamma}):=C_p(\epsilon) r_p({\gamma}_{-n}) e^{-\epsilon/2}$ and
 $\varphi =1$ yields a measurable function $\widehat{r}_p : \widehat{\mathcal Y} \to ]0,1]$ such that
 
 \begin{eqnarray}\label{R}
 C_p(\epsilon)r_p(\gamma_{-n}) e^{-\epsilon/2} \ge \widehat{r}_p (\widehat{\gamma}) e^{-n\epsilon/2};\;\;
 \forall \widehat{\gamma}\in \widehat{\mathcal Y}.
 \end{eqnarray}

Let us now consider the sequence of functions $(v_n)_n$ which are defined
on $\widehat{\mathcal Y}$ by setting $v_n(\widehat{\gamma}):=e^{-nL'}\prod_{j=1}^n e^{\widehat{u}_p(\gamma_{-j})}$.
By (\ref{Bir}), we have $\lim_n \frac{1}{n}\ln v_n(\widehat{\gamma}) =0$ for every $\widehat{\gamma}\in \widehat{\mathcal Y}$.
We may thus apply Lemma \ref{lemDD} to the functions $v_n$ and $\varphi =1$ and
get  a measurable function $\widehat{l}_p : \widehat{\mathcal Y} \to [1,+\infty[$ such
that $v_n\le \widehat{l}_p e^{n\epsilon/6}$, which means that
 
 \begin{eqnarray}\label{RR}
\prod_{j=1}^n \left(r_p(\gamma_{-j})\right)^{-1/2}= \prod_{j=1}^n e^{\widehat{u}_p(\gamma_{-j})}
\le \widehat{l}_p(\widehat{\gamma}) e^{nL'+n\epsilon/6};\;\; \forall \widehat{\gamma}\in \widehat{\mathcal Y}.
 \end{eqnarray}

We  now set $\widehat{\eta}_p := \frac{\widehat{r}_p}{\widehat{l}_p}$ and prove by induction on $n$ that:
\begin{enumerate}
\item $(f^p)_{\widehat{\gamma}}^{-n}$ is defined on
$T(\gamma_0, \widehat{\eta}_p(\widehat{\gamma})) \cap (U_0\times \Pj^k)$
 \item $\mbox{Lip} (f^p)_{\widehat{\gamma}}^{-n} \le \widehat{l}_p (\widehat{\gamma}) e^{n(L'+\tau + \epsilon/2)}$
 \item\label{item_incl} $(f^p)_{\widehat{\gamma}}^{-n}\left[T(\gamma_0, \widehat{\eta}_p(\widehat{\gamma}))
 \cap (U_0\times \Pj^k)\right] \subset T(\gamma_{-n},
C_p(\epsilon) r_p(\gamma_{-(n+1)})) \cap(U_0\times \Pj^k)$. 
\end{enumerate}

Let us stress that, setting $A=-(L'+\tau+\epsilon/2)$, the estimate on
$\mbox{Lip} (f^p)_{\widehat{\gamma}}^{-n}$ implies our statement  
since $\widehat{\eta}_p(\widehat{\gamma}) \widehat{l}_p(\widehat{\gamma})=\widehat{r}_p(\widehat{\gamma}) \le 1$.

Let us first check that this is true for $n=1$. By (\ref{R}) we have
$\widehat{\eta}_p(\widehat{\gamma})\le \widehat{r}_p (\widehat{\gamma}) \le C_p(\epsilon) r_p(\gamma_0)$
and therefore,
by Lemma \ref{Lemr},  $(f^p)_{\widehat{\gamma}}^{-1}$ is defined on
$T(\gamma_0, \widehat{\eta}_p(\widehat{\gamma})) \cap (U_0\times \Pj^k)$. Using Lemma
\ref{Lemr} again and (\ref{RR})
we get
$\mbox{Lip} (f^p)_{\widehat{\gamma}}^{-1}
\le e^{\tau +\epsilon/3} r_p(\gamma_{-1})^{\frac{-1}{2}}\le e^{L'+\tau + \epsilon/2} \widehat{l}_p (\widehat{\gamma}) $.
To show \ref{item_incl}
it suffices
to check that $ \widehat{l}_p (\widehat{\gamma})\widehat{\eta}_p(\widehat{\gamma})e^{L'+\tau + \epsilon/2}
= \widehat{r}_p (\widehat{\gamma})e^{L'+\tau + \epsilon/2}$
is smaller than $C_p(\epsilon)r_p(\gamma_{-2})$ which follows from (\ref{R})
since $e^{L'+\tau +\epsilon} \le 1$.

Let us now proceed with the inductive step.
Since $(f^p)_{\widehat{\gamma}}^{-(n+1)}= (f^p)_{\gamma_{-(n+1)}}^{-1)} \circ 
(f^p)_{\widehat{\gamma}}^{-(n)}$, Lemma \ref{Lemr} shows that
$(f^p)_{\widehat{\gamma}}^{-(n+1)}$ is well defined on $T(\gamma_0, \widehat{\eta}_p(\widehat{\gamma})) \cap (U_0\times \Pj^k)$.
To compute the Lipschitz constant we will use Lemma \ref{Lemr} and (\ref{RR}):
\begin{eqnarray*}
\mbox{Lip} (f^p)_{\widehat{\gamma}}^{-(n+1)} \le \prod_{j=1}^{n+1}
\mbox{Lip} (f^p)_{\gamma_{-j}}^{-1} = \prod_{j=1}^{n+1} e^{\tau +\epsilon/3} r_p(\gamma_{-j})^{-\frac{1}{2}}\\
\le e^{(n+1)(\tau +\epsilon/3)}\widehat{l}_p(\widehat{\gamma}) e^{(n+1)L'+(n+1)\epsilon/6}=\widehat{l}_p(\widehat{\gamma})
e^{(n+1)(L'+\tau +\epsilon/2)}.
\end{eqnarray*}
Finally,
\begin{eqnarray*}
\begin{aligned}
(f^p)_{\widehat{\gamma}}^{-{n+1}}\left[T(\gamma_0, \widehat{\eta}_p(\widehat{\gamma}))
\cap (U_0\times \Pj^k)\right]
&\subset T(\gamma_{-{n+1}}, \mbox{Lip}
(f^p)_{\widehat{\gamma}}^{-(n+1)} \widehat{\eta}_p(\widehat{\gamma})) \cap(U_0\times \Pj^k)\\
&\subset
T(\gamma_{-(n+1)},
C_p(\epsilon) r_p(\gamma_{-(n+2)} )) \cap (U_0\times \Pj^k),
\end{aligned}
\end{eqnarray*}
where the last inclusion follows from (\ref{R}), the estimate on
$\mbox{Lip} (f^p)_{\widehat{\gamma}}^{-(n+1)}$ and the fact that $e^{L'+\tau +\epsilon}\le 1$.
\fin

\subsection{Estimate of a Lyapunov exponent}

As Proposition \ref{PropInter} shows,  Theorem \ref{LemSF} will be proved if we
establish first that the functions $\widehat{u}_n$ are $\widehat{\mathcal M}$-integrable and then that
 $\lim_n\frac{1}{n} \int_{\widehat{X}} \widehat{u}_n\;d\widehat{\mathcal M}=L$ for
 some $L\le \frac{-\ln d}{2}$.
  This subsection is devoted to the proof of these two facts. We refer to Proposition
  \ref{ThmLya} below for a precise statement
  which, in particular, shows that the constant $L$ may be considered as a bound for a
  Lyapunov exponent of the system $({\mathcal J},{\mathcal F},{\mathcal M})$.
  
In the next Lemma, we list some basic properties of the functions $u_n$ and $\widehat{u}_n$.

\begin{lem}\label{LemPty} Let $U_0\Subset V_0\Subset W_0$ be open balls
centered at $\la_0$ in $M$. Let  $u_n : {\mathcal X}\times \overline{W_0} \to \R$
and $\widehat{u}_n : {\mathcal X}\to\R$ be defined by 
$u_n(\gamma,\la):=\ln \Vert  (DF_{\gamma(\la)}^n(0))^{-1}\Vert$ and
${\widehat u}_n(\widehat{\gamma}):=\sup_{\la \in U_0} u_n(\gamma_0,\la)$.
Let $\chi_1(\la)$ is the smallest Lyapunov exponent of the system $(J_\la,f_\la,\mu_\la)$.
The functions $u_n$ satisfy the following properties.
\begin{itemize}
\item[1)]  $u_n(\gamma,\cdot)$ is $p.s.h$ on $W_0$ for every $\gamma \in {\mathcal X}$.
\item[2)] The sequence $({\widehat u}_n)_n$ is subadditive i.e.
${\widehat u}_{m+n}\le {\widehat u}_n+ {\widehat u}_m\circ \widehat{\mathcal F}^n$.
\item[3)] For any fixed $\la \in W_0$, we have $\lim_n \frac{1}{n} u_n(\gamma,\la) =-\chi_1(\la)$
for ${\mathcal M}$-almost every $\gamma\in {\mathcal X}$.
\item[4)] For $\mathcal M$-almost every $\gamma\in {\mathcal X}$ we have
$\lim_n \frac{1}{n} u_n(\gamma,\la)=-\chi_1(\la)$ for Lebesgue-almost every $\la~\in~W_0$.
\end{itemize}
\end{lem}

\proof 1) When $\gamma\in {\mathcal X}$ is fixed the function $u_n(\gamma,\cdot)$ is clearly continuous 
on $W_0$ and $u_n(\gamma,\la)=\sup_{\Vert e \Vert =1} \ln \Vert  (DF_{\gamma(\la)}^n(0))^{-1}\cdot e\Vert$.
To see that 
$u_n(\gamma,\cdot)$ is $p.s.h$ it thus suffices
to check that the map $\la~\mapsto~\ln\Vert(DF_{\gamma(\la)}^n(0))^{-1}\cdot e\Vert$ is $p.s.h$
for each unit vector
$e\in \C^k$, which is clear since
$\ln \Vert  (DF_{\gamma(\la)}^n(0))^{-1}~\cdot~e\Vert= \frac{1}{2} \ln \sum_{i=1}^k \vert a_i(\la)\vert ^2$, where the functions
$a_i$ are holomorphic on $\overline {W_0}$.

2)
Let $\gamma \in {\mathcal X}$ and $m,n\ge 1$. By definition, for every
$\la \in {W_0}$ and $z$ close enough to the origin we have:

\begin{eqnarray*}
F_{{\mathcal F}^n(\gamma)(\la)}^m\circ F_{\gamma(\la)}^n =
\left(\psi_{i({\mathcal F}^{m+n}(\gamma)),{\mathcal F}^{m+n}(\gamma)(\la)}\right)^{-1} \circ f_\la^m
\circ \psi_{i({\mathcal F}^{n}(\gamma)),{\mathcal F}^{n}(\gamma)(\la)}\circ\\
\circ\left(\psi_{i({\mathcal F}^{n}(\gamma)),{\mathcal F}^{n}(\gamma)(\la)} \right)^{-1} \circ f_\la^n \circ
\psi_{i(\gamma),\gamma(\la)} (z) = F_{\gamma(\la)} ^{m+n} (z)
\end{eqnarray*}
and therefore  
\begin{eqnarray}\label{ChRu}
\left(DF_{\gamma(\la)} ^{m+n} (0)\right)^{-1} = \left(DF_{\gamma(\la)} ^{n} (0)\right)^{-1}
\circ \left(DF_{{\mathcal F}^n(\gamma)(\la)} ^{m} (0)\right)^{-1}.
\end{eqnarray}
Thus, if $\widehat{\gamma}\in \widehat{\mathcal X}$ we have
\begin{eqnarray*} 
\widehat{u}_{m+n}(\widehat{\gamma}) = \ln \sup_{\la\in {U_0}} \Vert (DF_{\gamma_0(\la)} ^{m+n} (0))^{-1} \Vert=
\ln \sup_{\la\in {U_0}} \Vert (DF_{\gamma_0(\la)} ^{n} (0))^{-1} \circ  (DF_{{\mathcal F}^n(\gamma_0)(\la)} ^{m} (0))^{-1}\Vert\\
\le
 \ln \sup_{\la\in {U_0}} (\Vert (DF_{\gamma_0(\la)} ^{n} (0))^{-1}\Vert \;\Vert  (DF_{{\mathcal F}^n(\gamma_0)(\la)} ^{m} (0))^{-1}\Vert)
\\ \le 
\ln  \sup_{\la\in {U_0}} \Vert (DF_{\gamma_0(\la)} ^{n} (0))^{-1}\Vert
+\ln \sup_{\la\in {U_0}} \Vert  (DF_{{\mathcal F}^n(\gamma_0)(\la)} ^{m} (0))^{-1}\Vert 
=  {\widehat u}_n(\widehat \gamma)+ {\widehat u}_m(\widehat{\mathcal F}^n(\widehat \gamma)).
\end{eqnarray*}

3)  By Oseledec theorem (see \cite{KH_intro}), the subset $J_{\la,1}$ of $J_\la\setminus C_{f_{\la}}$ defined by 
\begin{center}
$ J_{\la,1}:=\{x\in J_\la\setminus C_{f_\la}\;\colon\; \lim_n \frac{1}{n} \ln \Vert (Df_\la^n)_x^{-1}\Vert =-\chi_1(\la)\}$
\end{center}

has full $\mu_\la$ measure. As $p_{\la\star}\left({\mathcal M}\right)=\mu_\la$,
this implies that $\gamma(\la) \in J_{\la,1}$ for $\mathcal M$-almost every $\gamma$ in $\mathcal X$
and the assertion follows since
$\lim_n \frac{1}{n} \ln \Vert (Df_\la^n)_{\gamma(\la)}^{-1}\Vert= \lim_n \frac{1}{n}
\ln \Vert (DF_{\gamma(\la)}^n(0))^{-1}\Vert=\lim_n \frac{1}{n} u_n(\gamma,\la)$.

4) Let us denote by $\mathcal{L}$ the Lebesgue measure on $M$. Let $E$ be the
${\mathcal M}\otimes{\mathcal L}$-measurable subset of ${\mathcal X}\times W_0$ given by
$$E:=\{(\gamma,\la)\in {\mathcal X}\times W_0\;\colon\; \lim_n\frac{1}{n} u_n(\gamma,\la)=-\chi_1(\la)\}.$$
For every $\la\in W_0$ and every $\gamma\in {\mathcal X}$ we set
$$E^\la:=\{\gamma\in{\mathcal X}\;\colon\; (\gamma,\la)\in E\}\;\;\;\textrm{and}\;\; \;E_{\gamma}:=
\{\la\in W_0\;\colon\; (\gamma,\la)\in E\}.$$
We have to show that ${\mathcal L}(E_\gamma)={\mathcal L}(W_0)$ for $\mathcal M$-almost every $\gamma \in {\mathcal X}$.
Since, according to the above third assertion, ${\mathcal M}(E^\la)=1$
for every $\la\in W_0$, this immediately follows from Tonelli's theorem:
$$\int_{\mathcal X} {\mathcal L}(E_\gamma)\;d{\mathcal M}(\gamma)= {\mathcal M}\otimes{\mathcal L}(E)
=\int_{W_0} {\mathcal M}(E^\la)\;d{\mathcal L}(\la) = {\mathcal L} (W_0).$$\fin

Our strategy will be to transfer the  estimates known for the
system $(J_{\la_0},f_{\la_0},\mu_{\la_0})$ to the system $({\mathcal X},{\mathcal F},{\mathcal M})$.
This is  possible because the graphs $\Gamma_{\gamma}$ for $\gamma\in {\mathcal X}$
must approach the critical set $C_f$ locally uniformly, a phenomenon which simply
relies on the compactness of the closure of $\mathcal X$ and the following
basic property (see \cite{BD_stab}).\\

{\bf Fact} 
{\it There exist $0<\alpha\le 1$ such that $\sup_{V_0} \vert
\varphi\vert \le \vert \varphi(\la)\vert ^{\alpha}$ for every $\la\in V_0$ and every holomorphic function
$\varphi : W_0\to\C$ such that $0<\vert \varphi\vert<1$.\\}

More specifically, the key uniformity property we need is given by the next  lemma.
In our proofs, we shall denote the smallest singular value of an invertible linear map $L$ of $\C^k$ by $\delta(L)$. 
Let us recall  that $\delta(L)=\Vert L^{-1}\Vert^{-1}$ and that
$\vert \mbox{det} L \vert \ge \delta(L) ^k$.

\begin{lem}\label{LemKey}
Let $U_0$,$V_0$,$W_0$ and $u_n : {\mathcal X}\times \overline{W_0} \to \R$ be 
as in Lemma \ref{LemPty}. Then there exist $\alpha >0$ and $c>0$
such that $\frac{1}{n} u_n(\gamma,\la) \le \frac{k}{\alpha} \frac{1}{n}u_n(\gamma,\la'_0) +\ln c$
for every $n\geq 1$, every $\gamma \in {\mathcal X}$ and every $\la'_0, \la\in {V_0}$.
\end{lem}

\proof
From the compactness of $\overline{\mathcal X}$ and $\overline{V_0}$ we get the
existence of a constant $c_1 >0$ such that $\vert \mbox{det} (DF_{\gamma(\la)}^1(0))\vert \le c_1 \delta  (DF_{\gamma(\la)}^1(0))$ 
for every $\la\in \overline{V_0}$ and every $\gamma\in {\mathcal X}$.

Then, as $\mbox{det} DF_{\gamma(\la)}^n(0)=\prod _{j=0}^{n-1} \mbox{det} DF_{{\mathcal F}^j(\gamma)}^1(0)$  and 
$ \prod _{j=0}^{n-1} \delta( DF_{{\mathcal F}^j(\gamma)}^1(0)) \le  \delta ( DF_{\gamma(\la)}^n(0))$ 
we get 
\begin{eqnarray}\label{hh}
\vert \mbox{det} DF_{\gamma(\la)}^n(0)\vert \le c_1^n  \delta ( DF_{\gamma(\la)}^n(0));\;\;
\forall \gamma\in{\mathcal X},\;\forall \la\in\overline{V_0}.
\end{eqnarray}

Let us set $c_2:=\sup_{\la\in \overline{W_0},\gamma\in \overline{\mathcal X}} \vert \mbox{det} DF_{\gamma(\la)}^1 (0)\vert$.
When $\gamma \in {\mathcal X}$,  the holomorphic function 
$\varphi(\la):=\frac{1}{c_2^n} \mbox{det} DF_{\gamma(\la)}^n (0)$ is non vanishing and
its modulus is bounded by $1$ on $W_0$. Applying the above stated Fact  to $\varphi$, we get $0<\alpha\le 1$ 
(which only depends on $V_0$ and $W_0$) such that:
\begin{eqnarray}\label{h}
\sup_{\la\in V_0} \vert \mbox{det} DF_{\gamma(\la)}^n(0)\vert
\le c_2^{n(1-\alpha)} \vert \mbox{det} DF_{\gamma(\la)}^n(0)\vert^{\alpha};\;\;\forall n\ge 1,\;
\forall \gamma\in {\mathcal X},\;\forall \la\in V_0.
\end{eqnarray}

Using successively (\ref{h}) and (\ref{hh}) we get for any $\la,\la'_0\in V_0$
\begin{eqnarray*}
\begin{aligned}
\left[\delta(DF_{\gamma(\la'_0)}^n(0))\right]^k
\le \vert \mbox{det} DF_{\gamma(\la'_0)}^n(0))\vert
&\le c_2^{n(1-\alpha)} \vert \mbox{det} DF_{\gamma(\la)}^n(0)\vert^\alpha\\
& \le  c_2^{n(1-\alpha)} c_1^{n\alpha} \left[ \delta( DF_{\gamma(\la)}^n(0))\right]^\alpha.
\end{aligned}
\end{eqnarray*}

Then, applying $\ln$ and multiplying by $\frac{-1}{n}$ we get 
$$k\frac{1}{n}u_n(\gamma,\la'_0) \ge \alpha \frac{1}{n} u_n(\gamma,\la) -\alpha \left(\ln c_1 +\frac{1-\alpha}{\alpha} \ln c_2\right)$$
which is the desired estimate with 
$c:=c_1 c_2^{(1-\alpha)/\alpha} $.\fin

The next Lemma gathers the properties of the sequence $(u_n)_n$ which will be crucial to end our proof.

\begin{lem}\label{PropCru}Let $U_0$,$V_0$,$W_0$ and $u_n : {\mathcal X}\times \overline{W_0} \to \R$,
$\widehat{u}_n : \widehat{\mathcal X}\to\R$ be 
as in Lemma \ref{LemPty}. Then the following properties occur.
\begin{itemize}
\item[1)] The sequence $(\frac{1}{n} u_n)_n$ is uniformly bounded from below on ${\mathcal X}\times V_0$.
\item[2)] The sequence $\frac{1}{n} u_n(\gamma,\cdot)$ is uniformly bounded on $V_0$
for $\mathcal M$-almost every $\gamma\in {\mathcal X}$.
\item[3)] The functions ${\widehat u}_n$ are $\widehat{\mathcal M}$-integrable.
\end{itemize}
\end{lem}

\proof
1) Using the properties of the smallest singular value we have
\begin{eqnarray*}
\begin{aligned}
\frac{1}{n} u_n(\gamma,\la) = \frac{-1}{n} \ln \delta \left(DF_{\gamma(\la)}^n (0)\right)
& \ge \frac{-1}{nk} \ln \left| \mbox{det} \left(DF_{\gamma(\la)} ^n (0)\right)\right| \\
& =\frac{-1}{k} \left(\frac{1}{n} \sum_{j=0}^{n-1} \ln \left| \mbox{det} DF_{{\mathcal F}^j(\gamma)(\la)}^1 (0)\right|\right)
\end{aligned}
\end{eqnarray*}
and the assertion follows immediately from the definition and the continuity of $F_{\gamma(\la)}^1$.

2) 
We have just seen that $\frac{1}{n} u_n(\gamma,\cdot)$ is uniformly bounded from below on $V_0$.
By the fourth assertion of Lemma  \ref{LemPty}, for ${\mathcal M}$-almost every $\gamma\in{\mathcal X}$
there exists $\la_\gamma\in V_0$ such that
$\lim_n \frac{1}{n} u_n(\gamma,\la_\gamma)=-\chi_1(\la_\gamma)$. On the other
hand, by Lemma \ref{LemKey}, we have $\frac{1}{n} u_n(\gamma,\la)
\le \frac{k}{\alpha} \frac{1}{n} u_n(\gamma,\la_\gamma) +\ln c$ for every $n\in \N$ and
every $\la \in V_0$ and thus $\frac{1}{n} u_n(\gamma,\cdot)$ is uniformly bounded from above on $V_0$.

3) By the above first assertion, we know that $\widehat{u}_n$ is bounded from below. It thus suffices to show that 
$\int \widehat {u}_n(\widehat{\gamma}) \;d\widehat{\mathcal M} (\widehat{\gamma}) < +\infty$.  
By Lemma \ref{LemKey} we have
\begin{eqnarray*}
\int \widehat {u}_n(\widehat{\gamma}) \;d\widehat{\mathcal M} (\widehat{\gamma})
\le n\ln c + \frac{k}{\alpha} \int {u}_n((\pi_0({\widehat \gamma}),\la_0) \;d\widehat{\mathcal M} (\widehat{\gamma})
= n\ln c + \frac{k}{\alpha} \int {u}_n( \gamma,\la_0) \;d{\mathcal M} ({\gamma})\\
=n\ln c + \frac{k}{\alpha} \int \ln \Vert (DF_{\gamma(\la_0)}^n (0) )^{-1}\Vert \;d{\mathcal M} (\gamma)
= n\ln c - \frac{k}{\alpha} \int \ln \delta (DF_{\gamma(\la_0)}^n (0) ) \;d{\mathcal M} (\gamma).\\
\end{eqnarray*}

Using (\ref{hh}), we thus get

\begin{eqnarray*}
\int \widehat {u}_n(\widehat{\gamma}) \;d\widehat{\mathcal M} (\widehat{\gamma})
\le
 -\frac{k}{\alpha} \int \ln \vert \mbox{det} (DF_{\gamma(\la_0)}^n (0) )\vert
 \;d{\mathcal M} (\gamma) + \frac{kn}{\alpha} \ln c_1 + n\ln c\\
 =- \frac{k}{\alpha} \int \ln  \vert \mbox{det} (Df_{\la_0}^n)_{\gamma(\la_0)}\vert  \;\;d{\mathcal M} (\gamma)  + C 
 =- \frac{k}{\alpha} \int \ln  \vert \mbox{det} (Df_{\la_0}^n)_x \vert \;(dp_{\la_0 \star}{\mathcal M})(x) + C
\end{eqnarray*}
and the conclusion follows from the integrability of
$\ln  \vert \mbox{det} (Df_{\la_0}^n)_x\vert$ with respect to $p_{\la_0\star} {\mathcal M}=\mu_{\la_0}$
(see \cite{DS_dyn}). 
\fin

We are now ready to establish the main result of this section.

\begin{prop}\label{ThmLya}  There exists a constant $L\le \frac{-\ln d}{2}$ such that 
$$\lim_n \frac{1}{n}\int_{\widehat{\mathcal X}} \widehat{u}_n\;d\widehat{\mathcal M} =\lim_n
\frac{1}{n}\int_{\widehat{\mathcal X}}\; \sup_{\la\in U_0}\ln \Vert  (DF_{\gamma_0(\la)}^n(0))^{-1}\Vert\;d\widehat{\mathcal M}= L$$
$\textrm{and}\;\lim_n \widehat{u}_n(\widehat{\gamma})=\lim_n \sup_{\la\in U_0}
\ln \Vert  (DF_{\gamma_0(\la)}^n(0))^{-1}\Vert=L\;\textrm{for}\;\widehat {\mathcal M}\textrm{-almost every}\; 
\widehat{\gamma} \in {\mathcal X}.$
\end{prop}

\proof 
We will apply Kingman's subadditive ergodic theorem (see \cite{A_random,D_propr}) to the sequence
$(\widehat{u}_n)_n$. This is possible since the system $(\widehat {\mathcal X},\widehat {\mathcal F},\widehat {\mathcal M})$
is ergodic, the sequence $(\widehat{u}_n)_n$ is subadditive
(second assertion of Lemma \ref{LemPty}) and $\widehat{u}_1\in L^1(\widehat{\mathcal M})$
(last assertion of Lemma \ref{PropCru}). According to this theorem,
there exists $L\in \R$ such that $\lim_n \widehat{u}_n(\widehat{\gamma})=L$ for $\widehat{\mathcal M}$-almost every
 $\widehat\gamma\in \widehat{\mathcal X}$
and $\lim_n \frac{1}{n}\int_{\widehat{\mathcal X}} \widehat{u_n}\;d\widehat{\mathcal M} =L$.
It remains to show that $L\le \frac{-\ln d}{2}$.

Taking into account the fourth assertion of Lemma \ref{LemPty} and the
second assertion of Lemma \ref{PropCru}, we may  thus pick $\widehat{\gamma} \in \widehat{\chi}$ such that:

\begin{itemize}
\item[i)] $\lim_n \widehat{u}_n(\widehat{\gamma})=L$
\item[ii)] $\frac{1}{n} u_n(\gamma_0,\cdot)$ is uniformly bounded on $V_0$
\item[iii)] $\lim_n \frac{1}{n} u_n(\gamma_0,\la)=-\chi_1(\la)$ for Lebesgue-almost every $\la\in V_0$.
\end{itemize}

We proceed by absurd.  Assuming that $L > \frac{-\ln d}{2}$, we will reach
a contradiction with the fact that $\chi_1(\la)\ge \frac{\ln d}{2}$ for all $\la$ (see \cite{BD_exp}).
 Recalling that $\widehat{u}_n(\widehat{\gamma})= \sup_{\la \in U_0} u_n(\gamma_0,\la)$, 
 there exist $\la_{n_k}\in U_0$ and $\epsilon >0$ such that
 $\frac{1}{n_k} u_{n_k} (\gamma_0,\la_{n_k}) \ge \frac{-\ln d}{2} + \epsilon$.
 We may assume that $\la_{n_k}\to \la'_0\in \overline{U_0}$ and pick $r>0$ such
 that $B(\la_{n_k},r)\subset V_0$ for all $k\in \N$. Then, by the subharmonicity
 of $u_{n_k}(\gamma_0,\cdot)$ on $V_0$ (first assertion of Lemma \ref{LemPty}) we get:
 $$\frac{-\ln d}{2} + \epsilon \le \frac{u_{n_k}(\gamma_0,\la_{n_k})}{n_k}
 \le \frac{1}{\vert B(\la_{n_k},r)\vert} \int_{B(\la_{n_k},r)} \frac{u_{n_k}(\gamma_0,\la)}{n_k}$$
 which, by Lebesgue dominated convergence theorem, yields $$\frac{-\ln d}{2} + \epsilon
 \le \frac{1}{\vert B(\la'_{0},r)\vert} \int_{B(\la'_{0},r)}  -\chi_1(\la)$$
 and contradicts the fact that $\chi_1(\la)\ge \frac{\ln d}{2}$ for all $\la$. \fin

 \end{document}